\input amstex
\magnification=1200
\documentstyle{amsppt}
\NoRunningHeads
\NoBlackBoxes
\define\mYB{\operatorname{\widetilde mYB}}

\define\>{\right>}
\define\<{\left<}
\define\soa{\operatorname{\frak s\frak o}}
\topmatter
\title
Topics in hidden symmetries. VI.
\endtitle
\author Denis V.Juriev
\endauthor
\affil\eightpoint\it ``Thalassa Aitheria'' Research Center for
Mathematical Physics and Informatics,\linebreak
ul.Miklukho-Maklaya 20-180, Moscow 117437 Russia.\linebreak
E-mail: denis\@juriev.msk.ru\tenpoint
\endaffil
\date q-alg/9708028\enddate
\abstract In this note devoted to some aspects of the inverse problem
of rep\-re\-sen\-tation theory the attention is concentrated on
the interrelations between various algebraic structures (algebras with
operators) unraveled by different solutions of the same inverse problem.
The Lie algebras and Jordan triple systems with operators are considered.
\endabstract
\endtopmatter
\document
This article being a continuation of five previous parts [1] as
illustrates the general ideology presented in the review [2] as explicates
its new features. The at\-tention is now concentrated presumably on the
interrelations between various algebraic struc\-tures (algebras with
operators) unraveled by different solutions of the same inverse problem
of the representation theory.

\head 1. Topic Eleven: Jordan designs on Lie $\mYB$-algebras\endhead

\subhead 1.1. Lie $\mYB$-algebras and Jordan designs\endsubhead

\definition{Definition 1A {\rm [3]}} The {\it Lie $\mYB$-algebra\/} is a
pair $(\frak g,R)$, where $\frak g$ is the Lie algebra with brackets
$[\cdot,\cdot]$ and $R$ is an operator in it such that
$$R[RX,Y]+R[X,RY]=[RX,RY]+R^2[X,Y].$$
\enddefinition

The concept of the Lie $\mYB$-algebra was motivated by the constructions of
the paper [4], the main identity of the Lie $\mYB$-algebra coincides with the
modified classical Yang-Baxter equation if $R^2=1$ (such condition often holds
[4]). Examples of Lie $\mYB$-algebras were considered in [3].

\remark{Remark 1\/ {\rm [3]}} In any Lie $\mYB$-algebra $(\frak g,R)$ the bracket
$$[X,Y]_R=[RX,Y]+[X,RY]-R[X,Y]$$
supplies the linear space $\frak g$ by a structure of the Lie algebra.
\endremark

\remark{Remark 2\/ {\rm [3]}} Let $(\frak g,R)$ be a Lie $\mYB$-algebra and
$f(x)$ be a polynomial then $(\frak g,f(R))$ is the Lie $\mYB$-algebra, in
particular, $(\frak g,1+\lambda R)$ is a Lie $\mYB$-algebra for any $\lambda$.
Hence, the (compatible) Lie brackets $[\cdot,\cdot]_{1+\lambda R}$ supply the
linear space $\frak g$ by a structure of the linear bunch of Lie algebras.
\endremark

\definition{Definition 1B {\rm [5]}} The {\it Jordan design\/} on
the Lie algebra $\frak g$ with the commutator $[\cdot,\cdot]$ is
a $\frak g$-equivariant structure of a Jordan triple system on $\frak g$
with the trilinear operation $\<\cdot,\cdot,\cdot\>$ such that
$$[A,\<X,A,X\>]+[X,\<A,X,A\>]=0$$
for any $X$ and $Y$ from $\frak g$.
\enddefinition

Remind that the {\it Jordan triple system\/} [6] is the linear space $V$
supplied by a tri\-li\-near operation $\<\cdot,\cdot,\cdot\>$,
which obeys the following identity
$$\<x,\<a,z,b\>,y\>=\<\<x,a,y\>,b,z\>+\<\<y,a,z\>,b,x\>-\<\<x,b,y\>,a,z\>.$$

Examples of Jordan designs on Lie algebras were considered in [5]. Note
that both additional structures on Lie algebras (the structure of Lie
$\mYB$-algebra and the structure of Jordan design) naturally appear in the
context of the inverse problem of representation theory. So it is reasonable
to discuss their interrelations and to consider Jordan designs on
$\mYB$-Lie algebras now. It will be demonstrated how some nontrivial
conditions of the structural compatibility appear.

Let $(\frak g,R)$ be a Lie $\mYB$-algebra with the Jordan design
$\<\cdot,\cdot,\cdot\>$. Define a new triple product as
$$\aligned
\<X,Y,Z\>_R&=\<X,RY,RZ\>+\<RX,Y,RZ\>+\<RX,RY,Z\>\\
&-R\<RX,Y,Z\>-R\<X,RY,Z\>-R\<X,Y,RZ\>+R^2\<X,Y,Z\>.
\endaligned
$$
Note that $\<\cdot,\cdot,\cdot\>$ depends on $R$ quadratically whereas
$[\cdot,\cdot]_R$ depends on $R$ linearly. The construction of this new
product can be naturally generalized to the higher $n$-ary operations.

\definition{Definition 1C} The {\it Jordan triple $\mYB$-system\/} is a pair
$(\frak g,R)$, where $\frak g$ is the Jordan triple system with triple product
$\<\cdot,\cdot,\cdot\>$ and $R$ is an operator in it such that
$$R\<RX,Y,Z\>+R\<X,Y,RZ\>=\<RX,Y,RZ\>+R^2\<X,Y,Z\>.$$
\enddefinition

\remark{Remark 3} The operation $\<\cdot,\cdot,\cdot\>_R$ has the reduced
form
$$\<X,Y,Z\>_R=\<RX,RY,Z\>+\<X,RY,RZ\>-R\<X,RY,Z\>$$
in any Jordan triple $\mYB$-system.
\endremark

\proclaim{Theorem 1A} In any Jordan triple $\mYB$-system $(\frak g,R)$
the bracket $\<\cdot,\cdot,\cdot\>_R$ supp\-lies the linear space $\frak g$
by a new structure of the Jordan triple system.
\endproclaim

\remark{Remark 4} The equality
$$R\<X,Y,Z\>_R=\<RX,RY,RZ\>$$
holds in any Jordan triple $\mYB$-system.
\endremark

Let us now consider the Lie $\mYB$-algebra $(\frak g,R)$ supplied
by the $\frak g$-equivariant Jordan trilinear operation
$\<\cdot,\cdot,\cdot\>$ such that the pair $(\frak g,R)$ is a Jordan
triple $\mYB$-system. The following proposition holds.

\proclaim{Proposition 1} Let $(\frak g,R)$ be a Lie $\mYB$-algebra
and $\<\cdot,\cdot,\cdot\>$ be a $\frak g$-equivariant Jordan
trilinear operation on $\frak g$ such that $(\frak g,R)$ is a Jordan
$\mYB$-triple system. The Jordan trilinear operation
$\<\cdot,\cdot,\cdot\>_R$ is $\frak g_R$-equivariant, where $\frak g_R$
is the Lie algebra $\frak g$ supplied by the bracket $[\cdot,\cdot]_R$.
\endproclaim

\definition{Definition 1D} The {\it Jordan $\mYB$-design\/} on the Lie
$\mYB$-algebra $(\frak g,R)$ is a Jordan design on the Lie algebra $\frak g$
such that $(\frak g,R)$ is a Jordan triple $\mYB$-system.
\enddefinition

Let us now formulate the main theorem on Jordan $\mYB$-designs.

\proclaim{Theorem 1B} If $\<\cdot,\cdot,\cdot\>$ is a Jordan
$\mYB$-design on the Lie $\mYB$-algebra $(\frak g,R)$ then the
bracket $\<\cdot,\cdot,\cdot\>_R$ realizes a Jordan design on
the Lie algebra $\frak g_R$.
\endproclaim

\remark{Example 1} Let $\frak g$ be the Lie algebra $\soa(3,\Bbb C)$ of all
skew-symmetric $n\times n$ matrices and $\<\cdot,\cdot\>$ be the natural
bilinear symmetric form on it. Put $RX=\<X_0,X\>$ and $\<X,Y,Z\>=
\<X,Y\>Z+\<Y,Z\>X$ then $(\frak g,R)$ is a Lie $\mYB$-algebra with
a Jordan $\mYB$-design.
\endremark

Jordan $\mYB$-designs on Lie $\mYB$-algebras may be used for a description
of nonhamiltonian interactions governed by Jordan designs [5] in external
(magnetic-type) fields.

\subhead 1.2. Lie bi-$\mYB$-algebras and Jordan designs\endsubhead
Let us now discuss the Jordan designs on Lie bi-$\mYB$-algebras.

\definition{Definition 2A} The {\it Lie bi-$\mYB$-algebra\/} is a Lie algebra
$\frak g$ with bracket $[\cdot,\cdot]$ supplied by two commuting operators
$R_1$ and $R_2$ such that $(\frak g,R_1)$ and $(\frak g,R_2)$ are the Lie
$\mYB$-algebras with identical brackets $[\cdot,\cdot]_{R_1}$ and
$[\cdot,\cdot]_{R_2}$. A Lie bi-$\mYB$-algebra $(\frak g,R_1,R_2)$
is called {\it even-tempered} if the identities
$$\aligned
[R_1X,R_2Y]+[R_2X,R_1Y]-R_1R_2[X,Y]&=[R_1^2X,Y]+[X,R_1^2Y]-R_1^2[X,Y],\\
[R_1X,R_2Y]+[R_2X,R_1Y]-R_1R_2[X,Y]&=[R_2^2X,Y]+[X,R_2^2Y]-R_2^2[X,Y]
\endaligned
$$
hold.
\enddefinition

Examples of (even-tempered) bi-$\mYB$-algebras were considered in [3].
The main example is exposed below.

\remark{Example 2} Let $\frak A$ be an associative algebra then its
commutator algebra $\frak A_{[\cdot,\cdot]}$ is a Lie bi-$\mYB$-algebra, where
$R_1(X)=R_Q^r(X)=XQ$ or $R_2(X)=R_Q^l(X)=QX$ are the operators of the left or
right multiplication on the elements $Q$ of the associative algebra $\frak A$.
Moreover,
$$[X,Y]_{R_1}=[X,Y]_{R_2}=XQY-YQX.$$
\endremark

The Lie $\mYB$-algebra $(\frak g,R)$ is a Lie bi-$\mYB$-algebra with
$R_1=R_2=R$. Such Lie bi-$\mYB$-algebra is not even-tempered in general.

\remark{Remark 5\/ {\rm [3]}} A Lie $\mYB$-algebra $(\frak g,R)$ is the Lie
bi-$\mYB$-algebra if and only if there exists a derivative $\xi$ of the Lie
algebra $\frak g$ commuting with $R$ such that
$$[\xi X,\xi Y]=[SX,Y]+[X,SY]-S[X,Y],\qquad S=R\xi,$$
and $R_1=R$, $R_2=R+\xi$.
\endremark

Note that the operator $\xi$ is a derivative of both brackets $[\cdot,\cdot]$
and $[\cdot,\cdot]_R$.

\remark{Remark 6\/ {\rm [3]}} For any Lie bi-$\mYB$-algebra
$(\frak g,R_1,R_2)$ and any polynomial $f(x)$ the triple
$(\frak g,f(R_1),f(R_2))$ is a Lie bi-$\mYB$-algebra.
\endremark

\remark{Remark 7\/ {\rm [3]}} The identities in the even-tempered Lie
bi-$\mYB$-algebra may be rewritten in terms of $R$ and $\xi$ as
$$[RX,\xi Y]+[\xi X,RY]-R\xi[X,Y]=[R^2X,Y]-2[RX,RY]+[X,R^2Y].$$
\endremark

Note that if one puts $R_0=\tfrac12(R_1+R_2)$ then the bracket
$[\cdot,\cdot]_{R_0}$ defined as $[X,Y]_{R_0}=[R_0X,Y]+[X,R_0Y]-R_0[X,Y]$
coincides with brackets $[\cdot,\cdot]_{R_i}$ ($i=1,2$), however, the pair
$(\frak g,R_0)$ does not constitute a Lie bi-$\mYB$-algebra.

Let us now define Jordan triple $\mYB$-triple systems.

\definition{Definition 2B} The {\it Jordan triple bi-$\mYB$-system\/}
is a Jordan triple system $\frak g$ with trilinear operation
$\<\cdot,\cdot,\cdot\>$ supplied by two commuting operators $R_1$ and $R_2$
such that $(\frak g,R_1)$ and $(\frak g,R_2)$ are the Jordan triple
$\mYB$-systems with identical trilinear operations
$\<\cdot,\cdot,\cdot\>_{R_1}$ and $\<\cdot,\cdot,\cdot\>_{R_2}$.
A Jordan triple bi-$\mYB$-system $(\frak g,R_1,R_2)$
is called {\it normal\/} and {\it even-tempered\/} if the identities
$$\aligned
\<X,R_1R_2Y,Z\>&=\<R_1X,Y,R_2Z\>+\<R_2X,Y,R_1Z\>-R_1R_2\<X,Y,Z\>\\
&=\<R_1X,R_1Y,Z\>+\<X,R_1Y,R_1Z\>-R_1\<X,R_1Y,Z\>\\
&=\<R_2X,R_2Y,Z\>+\<X,R_2Y,R_2Z\>-R_2\<X,R_2Y,Z\>
\endaligned
$$
and
$$\aligned
\<R_1X,R_1R_2Y,R_2Z\>&+\<R_2X,R_1R_2Y,R_2Z\>-R_1R_2\<X,R_1R_2Y,Z\>\\
&=\<R_1^2X,R_1^2Y,Z\>+\<X,R_1^2Y,R_1^2Z\>-R_1^2\<X,R_1^2Y,Z\>\\
&=\<R_2^2X,R_2^2Y,Z\>+\<X,R_2^2Y,R_2^2Z\>-R_2^2\<X,R_2^2Y,Z\>,
\endaligned
$$
respectively, hold.
\enddefinition

\remark{Remark 8} The defining identities for normal Jordan triple
$\mYB$-systems may be concisely rewritten as
$$\<X,Y,Z\>_R=\<X,\varrho Y,Z\>,$$
where $\varrho=R_1R_2$.
\endremark

\remark{Exercise} Rewrite the definition of Jordan triple $\mYB$-systems
and its proper sub\-classes of normal and even-tempered systems in terms
of operators $R=R_1$ and $\xi=R_2-R_1$.
\endremark

Note that any Jordan triple $\mYB$-system $(\frak g,R)$ is a Jordan triple
bi-$\mYB$-system with $R_1=R_2=R$. Such Jordan triple bi-$\mYB$-system is
neither normal nor even-tempered in general.

\remark{Example 3} Let $\frak A$ be an associative algebra then the
associated Jordan triple system $\frak A_{\<\cdot,\cdot,\cdot\>}$ (where
$\<X,Y,Z\>=XYZ+ZYX$) is a Jordan triple bi-$\mYB$-system, where
$R_1(X)=R_Q^r(X)=XQ$ or $R_2(X)=R_Q^l(X)=QX$ are the operators of the left or
right multiplication on the elements $Q$ of the associative algebra $\frak A$.
Moreover,
$$\<X,Y,Z\>_{R_1}=\<X,Y,Z\>_{R_2}=XQYQZ-ZQYQX.$$
\endremark

\proclaim{Proposition 2} The Jordan triple bi-$\mYB$-system of example 3 is
normal and even-tempered.
\endproclaim

Analizing the example 3 one may write a huge number of new identities
besides ones for normal and even-tempered Jordan triple $\mYB$-systems,
which hold for systems of this example, the most interesting of which is
the {\it $\varrho$-identity\/}
$$\<\varrho X,Y,\varrho Z\>=\varrho\<X,\varrho Y,Z\>.$$
For the triple systems of example 3 the mapping $\varrho:X\mapsto\varrho X$
is well-known quadratic mapping $X\mapsto QXQ=\tfrac12\<Q,X,Q\>$ [7].

\remark{Remark 9} In view of the remark 8 $\varrho$-identity may be concisely
rewritten as
$$\varrho\<X,Y,Z\>_{\varrho}=\<\varrho X,Y,\varrho Z\>,$$
where $\<\cdot,\cdot,\cdot\>_{\varrho}$ is the renotation for
the coinciding brackets $\<\cdot,\cdot,\cdot\>_{R_i}$ ($i=1,2$).
\endremark

\definition{Definition 2C} The {\it Jordan bi-$\mYB$-design\/} on the Lie
bi-$\mYB$-algebra $(\frak g,R_1,R_2)$ is a Jordan design on the Lie algebra
$\frak g$ such that $(\frak g,R_1,R_2)$ is a Jordan triple bi-$\mYB$-system.
\enddefinition

The design will be called normal or even-tempered if the Jordan triple
$\mYB$-system is normal or even-tempered.

\proclaim{Proposition 3} The Jordan triple bi-$\mYB$-system of example 3
supplies the even-tempered Lie bi-$\mYB$-algebra of example 2 by a normal
even-tempered Jordan bi-$\mYB$-design with $\varrho$-identity.
\endproclaim

\head 2. Topic Twelve: Lie $R\varrho$-algebras\endhead

\definition{Definition 3} The {\it Lie $R\varrho$-algebra\/} is a triple
$(\frak g,R,\varrho)$, where $\frak g$ is the Lie algebra with the bracket
$[\cdot,\cdot]$ and $R$, $\varrho$ are two operators in it such that
the following two identities
$$
\aligned
\varrho[X,Y]_{\varrho}&=[\varrho X,\varrho Y],\\
R[X,Y]_{\varrho}+\varrho[X,Y]_R&=[RX,\varrho Y]+[\varrho X,RY].
\endaligned
$$
holds for all $X$ and $Y$ from $\frak g$. Here
$$\aligned
[X,Y]_R&=[RX,Y]+[X,RY]-R[X,Y],\\
[X,Y]_{\varrho}&=[\varrho X,Y]+[X,\varrho Y]-\varrho[X,Y]+[RX,RY]-R[X,Y]_R.
\endaligned
$$
A Lie $R\varrho$-algebra $(\frak g,R,\varrho)$ is called {\it regular\/}
if the identity
$$R[X,Y]_R=2([\varrho X,Y]+[X,\varrho Y])$$
holds.
\enddefinition

\proclaim{Proposition 4} Any even-tempered Lie bi-$\mYB$-algebra
$(\frak g,R_1,R_2)$ is the regular Lie $R\varrho$-algebra with $R=R_1+R_2$,
$\varrho=R_1R_2$.
\endproclaim

The inverse statement is not true. Let us expose an example of the Lie
$R\varrho$-algebra, which is not a Lie bi-$\mYB$-algebra.

\remark{Example 4} Let $\frak g$ be the Lie algebra $\soa(n,\Bbb C)$ of
all skew-symmetric $n\times n$ matrices and $Q$ be a symmetric $n\times n$
matrix. If one puts $RX=QX+XQ$ and $\varrho X=QXQ$ then $(\frak g,R,\varrho)$
will be a Lie $R\varrho$-algebra, however, there are no any operators
$R_1$, $R_2$ in $\frak g$ constituting a Lie bi-$\mYB$-algebra $(\frak
g,R_1,R_2)$ such that $R_1+R_2=R$, $R_1R_2=\varrho$, in general.
\endremark

\proclaim{Proposition 4} The bracket $[\cdot,\cdot]_{\varrho}$ obeys
the Jacobi identity in any Lie $R\varrho$-algebra $(\frak g,R,\varrho)$.
\endproclaim

Let us now explicate the relation between Lie $R\varrho$-algebras and
quadratic $\Gamma$-bunches of Lie algebras analogous to one between
Lie $\mYB$-algebras and linear $\Gamma$-bunches of Lie algebras [3].

\definition{Definition 4 {\rm [3]}} The {\it (one-parametric) bunch of Lie
algebras\/} is a linear space $\frak g$ supplied by a family of Lie brackets
$[\cdot,\cdot]_\lambda$ ($\lambda\in\Bbb R$). The Lie algebra from $\frak g$,
which is defined by the bracket $[\cdot,\cdot]_\lambda$, is denoted by
$\frak g_\lambda$. A bunch of Lie algebras $\frak g$ is called the
{\it $\Gamma$-bunch} if there are defined homomorphisms $R_\lambda$ from
the Lie algebras $\frak g_\lambda$ into the Lie algebra $\frak g_0$,
otherwords, the equality
$$R_\lambda[X,Y]_\lambda=[R_\lambda X,R_\lambda Y],$$
where $[\cdot,\cdot]=[\cdot,\cdot]_0$, holds for all $X$ and $Y$ from
$\frak g$.
\enddefinition

Below there will be considered only smooth bunches, i.e. the bunches,
which Lie brackets form a smooth family, and smooth $\Gamma$-bunches,
for which the family of homomorphisms $R_\lambda$ is also smooth.

\proclaim{Proposition 5 {\rm (cf.[3])}} Let $\frak g$ be a smooth
$\Gamma$-bunch of Lie algebras. Define the tangent bracket $[\cdot,\cdot]_R$
and quadratic bracket $[\cdot,\cdot]_{\varrho}$ as
$$[X,Y]_R=\left.\frac{d[X,Y]_\lambda}{d\lambda}\right|_{\lambda=0}
\quad\text{and}\quad
[X,Y]_{\varrho}=\frac12\left.\frac{d^2[X,Y]_\lambda}{d\lambda^2}
\right|_{\lambda=0}.
$$
Then
$$
\aligned
[X,Y]_R&=[RX,Y]+[X,RY]-R[X,Y],\\
[X,Y]_{\varrho}&=[\varrho X,Y]+[X,\varrho Y]-\varrho[X,Y]+[RX,RY]-R[X,Y]_R,
\endaligned
$$
where
$$R=\left.\frac{dR_\lambda}{d\lambda}\right|_{\lambda=0}
\quad\text{and}\quad
\varrho=\frac12\left.\frac{d^2R_\lambda}{d\lambda^2}\right|_{\lambda=0}.
$$
\endproclaim

A bunch of Lie algebras $\frak g$ is called quadratic if the Lie brackets
$[\cdot,\cdot]_\lambda$ from a quadratic family, otherwords,
$[\cdot,\cdot]_\lambda=[\cdot,\cdot]+\lambda[\cdot,\cdot]_R+
\lambda^2[\cdot,\cdot]_\varrho$. A $\Gamma$-bunch of of Lie algebras is
called linear if not only Lie brackets but also the homomorphisms $R_\lambda$
form a quadratic family, i.e. $R_\lambda=1+\lambda R+\lambda^2\varrho$.

\proclaim{Theorem 2} There is a natural one-to-one correspondence between
quadratic $\Gamma$-bunches of Lie algebras and Lie $R\varrho$-algebras.
\endproclaim

\remark{Question} What does the regularity condition in Lie
$R\varrho$-algebras mean in terms of quadratic $\Gamma$-bunches of Lie
algebras?
\endremark

\Refs
\roster
\item"[1]" Juriev D., Topics in hidden symmetries. I-V. E-prints:
hep-th/9405050, q-alg/9610026, q-alg/9611003, q-alg/9611019,
funct-an/9611003.
\item"[2]" Juriev D., An excursus into the inverse problem of representation
theory. Report RCMPI-95/04 (August 1995) [e-version: mp\_arc/96-477].
\item"[3]" Juriev D., Topics in isotopic pairs and their representations.
III. Bunches of Lie algebras and modified classical Yang-Baxter equation.
E-print: q-alg/9708027.
\item"[4]" Semenov-Tian-Shanskii M.A., What the classical $r$-matrix is.
Funkts.anal.i ego prilozh. 17(4) (1983) 17-33.
\item"[5]" Juriev D., Jordan designs on Lie algebras and nonhamiltonian
interaction of Hamiltonian systems. E-print: mp\_arc/97-431.
\item"[6]" Jacobson N., Structure and representations of Jordan algebras.
Amer.Math.Soc., Pro\-vi\-dence, RI, 1971.
\item"[7]" Koecher M., Jordan algebras and applications. 1962.
\endroster
\endRefs
\enddocument